\documentclass[10pt,english]{smfart}

\usepackage[T1]{fontenc}
\usepackage[english,francais]{babel}
\usepackage{latexsym,amscd}
\usepackage{amsmath,amsfonts,amssymb,mathrsfs}
\usepackage{enumerate,euscript}
\usepackage{amssymb,url,xspace,smfthm}
\input xy
\xyoption{all}

\newcommand{\BibTeX}{{\scshape Bib}\kern-.08em\TeX}
\newcommand{\T}{\S\kern .15em\relax }
\newcommand{\AMS}{$\mathcal{A}$\kern-.1667em\lower.5ex\hbox
        {$\mathcal{M}$}\kern-.125em$\mathcal{S}$}

\DeclareMathOperator{\Pic}{Pic}

\DeclareMathOperator{\rang}{rk}

\DeclareMathOperator{\Spec}{Spec}

\def\esssup_#1{\underset{#1}{\mathrm{ess\,sup}}}

\title{Differentiability of the arithmetic volume function}
\date{\today}
\author{Huayi Chen}

\address{Universit\'e Paris Diderot --- Paris 7,
Institut de math\'ematiques de Jussieu, case 247, 4 place Jussieu, 75252 Paris Cedex}
\email{chenhuayi@math.jussieu.fr}

\begin{document}
\def\smfbyname{}

\begin{abstract}
We introduce the positive intersection product in
Arakelov geometry and prove that the arithmetic volume
function is continuously differentiable. As
applications, we compute the distribution function of
the asymptotic measure of a Hermitian line bundle and
several other arithmetic invariants.
\end{abstract}

\begin{altabstract}
On introduit le produit d'intersection positive en g\'eom\'etrie d'Arakelov et on d\'emontre que la fonction volume arithm\'etique est continuement d\'erivable. Comme applications, on calcule la fonction de r\'epartition de la mesure de probabilit\'e asymptotique d'un fibr\'e inversible hermitien ainsi que quelques d'autres invariants arithm\'etiques.
\end{altabstract}
\maketitle

\tableofcontents
\section{Introduction}
Let $K$ be a number field, $\mathcal O_K$ be its
integer ring and $\pi:X\rightarrow\Spec\mathcal O_K$ be
an arithmetic variety of relative dimension $d$. Recall
that the {\it arithmetic volume} of a continuous
Hermitian line bundle $\overline L$ on $X$ is by
definition
\begin{equation}\label{Equ:arithmetic volume}\widehat{\mathrm{vol}}(\overline
L):=\limsup_{n\rightarrow\infty}
\frac{\widehat{h}^0(X,\overline L^{\otimes n
})}{n^{d+1}/(d+1)!},\end{equation} where
\begin{equation*}\widehat{h}^0(X,\overline L^{\otimes n})
=\log\#\{s\in\pi_*(L^{\otimes n})\mid
\forall\,\sigma:K\rightarrow\mathbb C,\;
\|s\|_{\sigma,\sup}\leqslant 1\}.\end{equation*} The
properties of the arithmetic volume
$\widehat{\mathrm{vol}}$ (see
\cite{Moriwaki07,Moriwaki08,Yuan07,Yuan08,Chen_bigness,Chen_Fujita})
are quite similar to the corresponding properties of
the classical volume function in algebraic geometry.
Recall that if $Y$ is a projective variety defined over
a field $k$ and if $L$ is a line bundle on $Y$, then
the volume of $L$ is defined as
\[\mathrm{vol}(L):=\limsup_{n\rightarrow\infty}
\frac{\rang_k H^0(Y,L^{\otimes n})}{n^{\dim Y}/(\dim Y
)!}.\]

In \cite{Bou_Fav_Mat06}, Boucksom, Favre and Jonsson
have been interested in the regularity of the geometric
volume function. They have actually proved that the
function $\mathrm{vol}(L)$ is continuously
differentiable on the big cone. The same result has also been independently obtained by Lazarsfeld and Mus\c{t}at\v{a} \cite{Lazarsfeld_Mustata08}, by using Okounkov bodies. Note that the geometric
volume function is not second order derivable in
general, as shown by the blow up of $\mathbb P^2$ at a
closed point, see \cite[2.2.46]{LazarsfeldI} for
details. In the differential of $\mathrm{vol}$ appears
the positive intersection product, initially defined in
\cite{Boucksom_Demailly_Paun_Peternell} in the
analytic-geometrical framework, and redefined
algebraically in \cite{Bou_Fav_Mat06}.

Inspired by \cite{Bou_Fav_Mat06}, we introduce an analogue of
the positive intersection product in Arakelov geometry
and prove that the arithmetic volume function
$\widehat{\mathrm{vol}}$ is continuously differentiable
on $\widehat{\mathrm{Pic}}(X)$. We shall establish the
following theorem:
\begin{theo}\label{Thm:main theorem}
Let $\overline L$ and $\overline M$ be two continuous
Hermitian line bundles on $X$. Assume that $\overline
L$ is big. Then
\[D_{\overline L}\widehat{\mathrm{vol}}(\overline M):=\lim_{n\rightarrow+\infty}\frac{\widehat{\mathrm{vol}}
(\overline L^{\otimes n}\otimes\overline
M)-\widehat{\mathrm{vol}}(\overline L^{\otimes
n})}{n^d}\] exists in $\mathbb R$, and the function
$D_{\overline L}\widehat{\mathrm{vol}}$ is additive on
$\widehat{\mathrm{Pic}}(X)$. Furthermore, one has
\[D_{\overline L}\widehat{\mathrm{vol}}(\overline M)=
(d+1)\big\langle\widehat{c}_1(\overline L
)^d\big\rangle\cdot\widehat{c}_1(\overline M).\]
\end{theo}
Here the {\it positive intersection product}
$\big\langle\widehat{c}_1(\overline L )^d\big\rangle$
is defined as the least upper bound of self
intersections of ample Hermitian line bundles dominated
by $\overline L$ (see \S\ref{SubSec:Posit} {\it
infra}). In particular, one has
$\big\langle\widehat{c}_1(\overline L
)^d\big\rangle\cdot\widehat{c}_1(\overline
L)=\big\langle\widehat{c}_1(\overline L
)^{d+1}\big\rangle=\widehat{\mathrm{vol}}(\overline
L)$, which shows that the arithmetic Fujita
approximation is asymptotically orthogonal.

As an application, we calculate explicitly the
distribution function of the asymptotic measure (see
\cite{Chen08,Chen_bigness}) of a generically big
Hermitian line bundle in terms of positive intersection
numbers. Let $\overline L$ be a Hermitian line bundle
on $X$ such that $L_K$ is big. The asymptotic measure
$\nu_{\overline L}$ is the vague limit (when $n$ goes
to infinity) of Borel probability measures whose
distribution functions are determined by the filtration
of $H^0(X_K,L_K^{\otimes n})$ by successive minima (see
\eqref{Equ:nu} {\it infra}). Several asymptotic
invariants can be obtained by integration with respect
to $\nu_{\overline L}$. Therefore, it is interesting to
determine completely the distribution of
$\nu_{\overline L}$, which will be given in Proposition
\ref{Pro:distribution function} by using the positive
intersection product.

The article is organized as follows. In the second
section, we recall some positivity conditions for
Hermitian line bundles and discuss their properties. In
the third section, we define the positive intersection
product in Arakelov geometry. It is in the fourth
section that we establish the differentiability of the
arithmetic volume function. Finally in the fifth
section, we present applications on the asymptotic
measure and we compare our result to some known results
on the differentiability of arithmetic invariants.

\bigskip

\noindent{\bf Acknowledgement:} I would like to thank
R. Berman, D. Bertrand, J.-B. Bost, S. Boucksom, C.
Favre and V. Maillot for interesting and helpful
discussions. I am also grateful to M. Jonsson for remarks.

\section{Notation and preliminaries}

In this article, we fix a number field $K$ and denote
by $\mathcal O_K$ its integer ring. Let $\overline K$
be an algebraic closure of $K$. Let
$\pi:X\rightarrow\Spec\mathcal O_K$ be a projective and
flat morphism and $d$ be the relative dimension of
$\pi$. Denote by $\widehat{\mathrm{Pic}}(X)$ the group
of isomorphism classes of (continuous) Hermitian line
bundles on $X$. If $\overline L$ is a Hermitian line
bundle on $X$, we denote by $\pi_*(\overline L)$ the
$\mathcal O_K$-module $\pi_*(L)$ equipped with sup
norms.

In the following, we recall several notions about
Hermitian line bundles. The references are
\cite{Gillet-Soule,Zhang95,BGS94,Moriwaki00}.

Assume that $x\in X(\overline K)$ is an algebraic point
of $X$. Denote by $K_x$ the field of definition
of $x$ and by $\mathcal O_x$ its integer ring. The morphism $x:\Spec\overline K\rightarrow X$ gives rise to a point $P_x$ of $X$ valued in
$\mathcal O_x$. The
pull-back of $\overline L$ by $P_x$ is a Hermitian line
bundle on $\Spec\mathcal O_x$. We denote by
$h_{\overline L}(x)$ its normalized Arakelov degree,
called the {\it height} of $x$. Note that the height
function is additive with respect to $\overline L$.

Let $\overline L$ be a Hermitian line bundle on $X$. We
say that a section $s\in\pi_*(L)$ is {\it effective}
(resp. {\it strictly effective}) if for any
$\sigma:K\rightarrow\mathbb C$, one has
$\|s\|_{\sigma,\sup}\leqslant 1$ (resp.
$\|s\|_{\sigma,\sup}< 1$). We say that the Hermitian
line bundle $\overline L$ is {\it effective} if it
admits a non-zero effective section.

Let $\overline L_1$ and $\overline L_2$ be two
Hermitian line bundles on $X$. We say that $\overline
L_1$ is {\it smaller} than $\overline L_2$ and we
denote by $\overline L_1\leqslant\overline L_2$ if the
Hermitian line bundle $\overline
L_1^\vee\otimes\overline L_2$ is effective.

We say that a Hermitian line bundle $\overline A$ is
ample if $A$ is {\it ample}, $c_1(\overline A)$ is
semi-positive in the sense of current on $X(\mathbb C)$
and $\widehat{c}_1(\overline L|_Y)^{\dim Y}>0$ for any
integral sub-scheme $Y$ of $X$ which is flat over
$\Spec\mathcal O_K$. Here the intersection number
$\widehat{c}_1(\overline L|_Y)^{\dim Y}$ is defined in
the sense of \cite{Zhang95} (see Lemma 6.5 {\it loc.
cit.}, see also \cite{Zhang95b}). Note that there always exists an ample
Hermitian line bundle on $X$. In fact, since $X$ is
projective, it can be embedded in a projective space
$\mathbb P^N$. Then the restriction of $\mathcal
O_{\mathbb P^N}(1)$ with Fubini-Study metrics on $X$ is
ample. Note that the Hermitian line bundle $\overline
L$ thus constructed has strictly positive smooth
metrics. Thus, if $\overline M$ is an arbitrary
Hermitian line bundle with smooth metrics on $X$, then
for sufficiently large $n$, $\overline
M\otimes\overline L^{\otimes n} $ is still ample.

We say that a Hermitian line bundle $\overline N$ is
{\it vertically  nef} if the restriction of $N$ on each
fiber of $\pi$ is nef and $c_1(\overline N)$ is
semi-positive in the sense of current on $X(\mathbb
C)$. We say that $\overline N$ is {\it nef} if it is
vertically nef and $\widehat{c}_1(\overline N|_Y)^{\dim
Y}\geqslant 0$ for any integral sub-scheme $Y$ of $X$
which is flat over $\Spec\mathcal O_K$. By definition,
an ample Hermitian line bundle is always nef.
Furthermore, if $\overline A$ is an ample Hermitian
line bundle and if $\overline N$ is a Hermitian line
bundle such that $\overline N^{\otimes
n}\otimes\overline A$ is ample for any integer
$n\geqslant 1$, then $\overline N$ is nef. We denote by
$\widehat{\mathrm{Nef}}(X)$ the subgroup of
$\widehat{\mathrm{Pic}}(X)$ consisting of nef Hermitian
line bundles.

If $f:X(\mathbb C)\rightarrow\mathbb R$ is a continuous
function, we denote by $\overline{\mathcal O}(f)$ the
Hermitian line bundle on $X$ whose underlying line
bundle is trivial, and such that the norm of the unit
section $\mathbf{1}$ at $x\in X(\mathbb C)$ is
$e^{-f(x)}$. Note that, if $f$ is positive, then
$\overline{\mathcal O}(f)$ is effective. If $f$ is
positive and plurisubharmonic, then $\overline{\mathcal
O}(f)$ is nef. In particular, for any $a\in\mathbb R$,
$\overline{\mathcal O}(a)$ is nef if and only if
$a\geqslant 0$. If $\overline L$ is a Hermitian line
bundle on $X$, we shall use the notation $\overline
L(f)$ to denote $\overline L\otimes\overline{\mathcal
O}(f)$.

We say that a Hermitian line bundle $\overline L$ is
{\it big} if its arithmetic volume
$\widehat{\mathrm{vol}}(\overline L)$ is strictly
positive. By \cite{Moriwaki07,Yuan07}, $\overline L$ is
big if and only if a positive tensor power of
$\overline L$ can be written as the tensor product of
an ample Hermitian line bundle with an effective one.
Furthermore, the analogue of Fujita's approximation
holds for big Hermitian line bundles, cf.
\cite{Chen_Fujita,Yuan08}.

The arithmetic volume function $\widehat{\mathrm{vol}}$
is actually a limit (cf. \cite{Chen_bigness}): one has
\[\widehat{\mathrm{vol}}(L)=\lim_{n\rightarrow\infty}
\frac{\widehat{h}^0(X,\overline L^{\otimes n
})}{n^{d+1}/(d+1)!}.\] Moreover, it is a birational
invariant which is continuous on
$\widehat{\mathrm{Pic}}(X)_{\mathbb Q}$, and can be
continuously extended to
$\widehat{\mathrm{Pic}}(X)_{\mathbb R}$, cf.
\cite{Moriwaki07,Moriwaki08}. The analogue of Siu's
inequality and the log-concavity hold for
$\widehat{\mathrm{vol}}$, cf. \cite{Yuan07,Yuan08}.

\begin{rema}
\begin{enumerate}[1)]
\item In \cite{Zhang95} and
\cite{Moriwaki00}, the notions of ample or nef line
bundles were reserved for line bundles with smooth
metrics, which is not the case here.
\item Note that there exists another (non-equivalent)
definition of arithmetic volume function in the
literature. See \cite[\S 10.1]{Ber_Bou08} and \cite[\S
5]{Chambert-Loir_Thuillier} where the ``arithmetic
volume'' of a Hermitian line bundle $\overline L$ was
defined as the following number:
\begin{equation}\label{Equ:sectional capa}S(\overline L):=\lim_{n\rightarrow+\infty}\frac{\chi(\pi_*(\overline
L^{\otimes n}))}{
n^{d+1}/(d+1)!}\in[-\infty,+\infty[,\end{equation}
which is also called {\it sectional capacity} in the
terminology of \cite{Rumely_Lau_Varley}. However, in
the analogy between Arakelov geometry and relative
algebraic geometry over a regular curve, it is
\eqref{Equ:arithmetic volume} that corresponds to the
geometric volume function. Note that one always has
\[\widehat{\mathrm{vol}}(\overline L)\geqslant \lim_{n\rightarrow+\infty}\frac{\chi(\pi_*(\overline
L^{\otimes n}))}{ n^{d+1}/(d+1)!},\] and the equality
holds when $\overline L$ is nef. Under this assumption, both quantities
are equal to the intersection number
$\widehat{c}_1(\overline L)^{d+1}$. This is a consequence of the Hilbert-Samuel formula. See \cite{Gillet-Soule,Abbes-Bouche,Zhang95,Autissier01,Randriam06,Moriwaki07} for details.
\end{enumerate}
\end{rema}

In the following, we present some properties of nef
line bundles. Note that Propositions \ref{Pro:critere
de nef} and \ref{Pro:positivity of nef intersection
produc} have been proved in \cite[\S 2]{Moriwaki00} for
Hermitian line bundles with smooth metrics. Here we
adapt these results to continuous metric case by using
the continuity of intersection numbers.
\begin{prop}\label{Pro:critere de nef}
Let $\overline N$ be a Hermitian line bundle on $X$
which is vertically nef. Assume that for any $x\in
X(\overline K)$, one has $h_{\overline N}(x)\geqslant
0$, then the Hermitian line bundle $\overline N$ is
nef.
\end{prop}
\begin{proof}
Choose an ample Hermitian line bundle $\overline A$ on
$X$ such that $h_{\overline A}$ has strictly positive
lower bound. For any integer $n\geqslant 1$, let
$\overline
L_n:=(L_n,(\|\cdot\|_\sigma)_{\sigma:K\rightarrow\mathbb
C })$ be the tensor product $\overline N^{\otimes
n}\otimes\overline A$. The height function
$h_{\overline L_n}$ is bounded from below by a strictly
positive number $\varepsilon_n$. Note that the metrics
of $\overline L_n$ are semi-positive. By \cite[Theorem
4.6.1]{Maillot00} (see also \cite[\S 3.9]{Randriam06}),
there exists a sequence of smooth positive metric
families $(\alpha_m)_{m\geqslant 1}$ with
$\alpha_m=(\|\cdot\|_{\sigma,m})_{\sigma:K\rightarrow\mathbb
C }$, such that $\|\cdot\|_{\sigma,m}$ converges
uniformly to $\|\cdot\|_{\sigma}$ when $m$ tends to the
infinity. Denote by $\overline L_{n,m}=(L_n,\alpha_m)$.
For sufficiently large $m$, $h_{\overline L_{n,m}}$ is
bounded from below by $\varepsilon_n/2$. Thus
\cite[Corollary 5.7]{Zhang95} implies that, for any
integer subscheme $Y$ of $X$ which is flat over
$\Spec\mathcal O_K$, one has $n^{-\dim Y
}\widehat{c}_1(\overline L_{n,m}|_Y )^{\dim Y}\geqslant
0$. By passing successively $m$ and $n$ to the
infinity, one obtains $\widehat{c}_1(\overline
N|_Y)^{\dim Y}\geqslant 0$. Therefore $\overline N$ is
nef.
\end{proof}

We say that a Hermitian line bundle $\overline L$ on
$X$ is {\it integrable} if there exist two ample
Hermitian line bundles $\overline  A_1$ and $\overline
A_2$ such that $\overline L=\overline
A_1\otimes\overline A_2^\vee$. Denote by
$\widehat{\mathrm{Int}}(X)$ the subgroup of
$\widehat{\mathrm{Pic}}(X)$ formed by all integrable
Hermitian line bundles. If $(\overline L_i)_{i=0}^d$ is
a family of integrable Hermitian line bundles on $X$,
then the intersection number
\[\widehat{c}_1(\overline L_0)\cdots\widehat{c}_1(\overline L_d)\]
is defined (see \cite[Lemma 6.5]{Zhang95}, \cite[\S
1]{Zhang95b} and \cite{Maillot00} \S 5). Furthermore,
it is a symmetric multi-linear form which is continuous
in each $\overline L_i $. Namely, for any family
$(\overline M_i)_{i=0}^d$ of integrable Hermitian line
bundles, one has
\[\lim_{n\rightarrow+\infty}n^{-d-1}\widehat{c}_1(
\overline L_0^{\otimes n}\otimes\overline
M_0)\cdots\widehat{c}_1(\overline L_d^{\otimes
n}\otimes\overline M_d)=\widehat{c}_1(\overline
L_0)\cdots\widehat{c}_1(\overline L_d).\]

\begin{prop}\label{Pro:positivity of nef intersection produc}
Let $(\overline L_i)_{i=0}^{d-1}$ be a family of nef
Hermitian line bundles on $X$ and $\overline M$ be an integrable Hermitian line bundle on $X$ which is effective. Then
\begin{equation}\label{Equ:positivity de intese}\widehat{c}_1(\overline
L_0)\cdots\widehat{c}_1(\overline L_{d-1})
\widehat{c}_1(\overline M)\geqslant 0\end{equation}
\end{prop}
\begin{proof}
Choose an ample Hermitian line bundle $\overline A$ on
$X$ such that $h_{\overline A}$ is bounded from below
by some strictly positive number. By virtue of the
proof of Proposition \ref{Pro:critere de nef}, for any
$i\in\{0,\cdots,d-1\}$ and any integer $n\geqslant 1$,
there exists a sequence of nef Hermitian line bundles
with smooth metrics $(\overline
L_{i,n}^{(m)})_{m\geqslant 1}$ whose underlying line
bundle is $L_i^{\otimes n}\otimes A$ and whose metrics
converge uniformly to that of $\overline L_i^{\otimes n
}\otimes\overline A$. By \cite[Proposition
2.3]{Moriwaki00}, one has \[\widehat{c}_1(\overline
L_0^{\otimes n}\otimes\overline
A)\cdots\widehat{c}_1(\overline L_{d-1}^{\otimes
n}\otimes \overline A)\widehat{c}_1(\overline
M)\geqslant 0.\] By passing to limit, one obtains
\eqref{Equ:positivity de intese}.
\end{proof}

\begin{rema}
Using the same method, we can prove that, if
$(\overline L_i)_{i=0}^d$ is a family of nef Hermitian
line bundles on $X$, then
\begin{equation}
\widehat{c}_1(\overline
L_0)\cdots\widehat{c}_1(\overline L_d)\geqslant 0.
\end{equation}
\end{rema}

\begin{prop}\label{Pro:critere de nefness}
Let $\overline{L}$ be a Hermitian line bundle on $X$
such that $c_1(\overline L)$ is semi-positive in the
sense of current on $X(\mathbb C)$. Assume that there
exists an integer $n> 0$ such that $L^{\otimes n}$ is
generated by its effective sections. Then the Hermitian
line bundle $\overline{L}$ is nef.
\end{prop}
\begin{proof}
Since $L^{\otimes n}$ is generated by its sections, the
line bundle $L$ is nef relatively to $\pi$. After
Proposition \ref{Pro:critere de nef}, it suffices to
verify that, for any $x\in X(\overline K)$, one has
$h_{\overline{L}}(x)\geqslant 0$. For any integer
$m\geqslant 1$, let $B_m=\pi_*(L^{\otimes m})$ and let
$B_m^{[0]}$ be the saturated sub-$\mathcal O_K$-module
of $B_m$ generated by effective sections. Since
$L^{\otimes n}$ is generated by its effective sections,
also is $L^{\otimes np}$ for any integer $p\geqslant
1$. In particular, one has surjective homomorphisms
$x^*\pi^* B_{pn,K_x}^{[0]}\rightarrow x^*
L_{K_x}^{\otimes np}$. By slope inequality (see
\cite[Appendix A]{BostBour96}), one has $nph_{\overline{
L}}(x)\geqslant\widehat{\mu}_{\min}(\overline{ B}_{np}
^{[0]}) $. By passing to limit, one obtains
$h_{\overline{L}}(x)\geqslant 0$.
\end{proof}

We say that a Hermitian line bundle $\overline L$ on
$X$ is {\it free} if $c_1(\overline L)$ is
semi-positive in the sense of current on $X(\mathbb C)$
and if some positive tensor power of $L$ is generated
by effective global sections. We denote by
$\widehat{\mathrm{Fr}}(X)$ the subgroup of
$\widehat{\mathrm{Pic}}(X)$ consisting of free
Hermitian line bundles. After Proposition
\ref{Pro:critere de nefness}, one has
$\widehat{\mathrm{Fr}}(X)\subset\widehat{\mathrm{Nef}}(X)$.

Unlike the ampleness, the properties of being big, nef,
or free are all invariant by birational modifications.
That is, if $\nu:X'\rightarrow X$ is a birational
projective morphism, and if $\overline L$ is a
Hermitian line bundle on $X$ which is big (resp. nef,
free), then also is $\nu^*(\overline L)$.

\section{Positive intersection product}

In this section, we shall define the positive
intersection product for big (non-necessarily integral)
Hermitian line bundles. When all Hermitian line bundles
are nef, the positive intersection product coincides
with the usual intersection product. Furthermore, the
highest positive auto-intersection number is just the
arithmetic volume of the Hermitian line bundle. We
shall use the positive intersection product to
interpret the differential of the arithmetic volume
function.

\subsection{Admissible decompositions}
\begin{defi}
Let $\overline L$ be a big Hermitian line bundle on
$X$. We call {\it admissible decomposition} of
$\overline L$ any triplet $(\nu,\overline N,p)$, where
\begin{enumerate}[1)]
\item $\nu:X'\rightarrow X$ is a birational projective
morphism,
\item $\overline N$ is a free Hermitian line bundle on $X'$,
\item $p\geqslant 1$ is an integer such that $\nu^*(\overline L^{\otimes p})
\otimes\overline N^\vee$ is effective.
\end{enumerate}
Denote by $\Theta(\overline L)$ the set of all
admissible decompositions of $\overline L$.
\end{defi}

We introduce an order relation on the set
$\Theta(\overline L)$. Let $D_i=(\nu_i:X_i\rightarrow
X,\overline N_i,p_i )$ ($i=1,2$) be two admissible
decompositions of $\overline L$. We say that $D_1$ is
{\it superior} to $D_2$ and we denote by $D_1\succ D_2$
if $p_2$ divides $p_1$ and if there exists a projective
birational morphism $\eta:X_1\rightarrow X_2$ such that
$\nu_2\eta=\nu_1$ and that $\overline
N_1\otimes(\eta^*\overline N_2)^{\vee\otimes(p_1/p_2)}$
is effective.

\begin{rema}\label{Rem:pull back of admissible decomp}
\begin{enumerate}[(1)]
\item Assume that $D=(\nu:X'\rightarrow X,
\overline{N},p)$ is an admissible decomposition of
$\overline{L}$. Then for any birational projective
morphism $\eta: X''\rightarrow X'$, the triplet
$\eta^*D:=(\nu\eta,\eta^*\overline{N},p)$ is also an
admissible decomposition of $\overline L$, and one has
$\eta^*D \succ D$.
\item Assume that $D=(\nu,
\overline{N},p)$ is an admissible decomposition of
$\overline{L}$. Then for any integer $n\geqslant 1$,
$D_n=(\nu,\overline{N}^{\otimes n}, np)$ is also an
admissible decomposition of $\overline{L}$.
Furthermore, one has $D_n\succ D$.
\item Assume that $D_1=(\nu,\overline{N}_1,p)$
and $D_2=(\nu,\overline{N}_2,q)$ are two admissible
decompositions of $\overline{L}$ whose underlying
birational projective morphisms are the same. Then
$D_1\otimes
D_2:=(\nu,\overline{N}_1\otimes\overline{N}_2,p+q)$ is
an admissible decomposition of $\overline{L}$.
\item Assume that $\overline M$ is an
effective Hermitian line bundle on $X$. By definition,
any admissible decomposition of $\overline L$ is also
an admissible decomposition of $\overline
L\otimes\overline M $.
\end{enumerate}
\end{rema}

In the following proposition, we show that the set
$\Theta(\overline L)$ is filtered with respect to the
order $\succ$.

\begin{prop}
\label{Pro:filtered} if $D_1$ and $D_2$ are two
admissible decompositions of $\overline L$, then there
exists an admissible decomposition $D$ of
$\overline{L}$ such that $D\succ D_1$ and $D\succ D_2$.
\end{prop}
\begin{proof}
After Remark \ref{Rem:pull back of admissible decomp}
(1)(2), we may assume that the first and the third
components of $D_1$ and $D_2$ are the same. Assume that
$D_1=(\nu,\overline{N}_1,p)$ and
$D_2=(\nu,\overline{N}_2,p)$, where $\nu: X'\rightarrow
X$ is a birational projective morphism. Let
$\overline{M}_i=\nu^*\overline{L}^{\otimes
p}\otimes\overline{N}_i^\vee$ ($i=1,2$). Since $\overline M_1$ and $\overline M_2$ are effective,
there exist homomorphisms $u_i:
M_i^\vee\rightarrow\mathcal O_{X'}$ corresponding to
effective sections $s_i:\mathcal O_{ X'}\rightarrow
M_i$ ($i=1,2$). Let $\eta:X''\rightarrow X'$ be the
blow up of the ideal sheaf $\mathrm{Im}(u_1\oplus u_2
)$. Let $ M$ be the exceptional line bundle and
$s:\mathcal O_{X''}\rightarrow M$ be the section which
trivializes $M$ outside the exceptional divisor. The
canonical surjective homomorphism $\eta^*(
M_1^\vee\oplus M_2^\vee )\rightarrow M^\vee$ induces by
duality an injective homomorphism $\varphi:
M\rightarrow M_1\oplus M_2$. We equip $  M_1\oplus M_2$
with metrics
$(\|\cdot\|_{\sigma})_{\sigma:K\rightarrow\mathbb C}$
such that, for any $x\in  X''_\sigma(\mathbb C)$ and
any section $(u,v)$ of $ M_{1,\sigma}\oplus
M_{2,\sigma} $ over a neighbourhood of $x$, one has
$\|(u,v)\|_{\sigma}(x)=\max\{\|u\|_{\sigma,1}(x),\|v\|_{\sigma,2}(x)\}$.
As $\varphi s=(\eta^*s_1,\eta^*s_2)$, and the sections
$s_1$ and $s_2$ are effective, one obtains that the
section $s$ is also effective. Let
$\overline{N}=(\nu\eta)^*\overline{L}^{\otimes p
}\otimes\overline{M}^\vee$. One has a natural
surjective homomorphism
\[\psi:\eta^*{N}_1\oplus
\eta^*{N}_2\longrightarrow N.\] Furthermore, if we
equip $\eta^*{N}_1\oplus \eta^*{N}_2$ with metrics
$(\|\cdot\|_{\sigma})_{\sigma:K\rightarrow\mathbb C}$
such that, for any $x\in X''_\sigma(x)$,
$\|(u,v)\|_\sigma(x)=\|u\|_\sigma(x)+\|v\|_{\sigma}(x)$,
then the metrics on $N$ are just the quotient metrics
by the surjective homomorphism $\psi$, which are
semi-positive since the metrics of $\eta^*\overline
N_1$ and of $\eta^*\overline N_2$ are. As both
Hermitian line bundles $N_1$ and $N_2$ are generated by
effective global sections, also is $\overline{N}$.
Therefore, $(\nu\eta,\overline{N},p)$ is an admissible
decomposition of $\overline{L}$, which is superior to
both $D_1$ and $D_2$.
\end{proof}

\subsection{Intersection of admissible decompositions}
Let $(\overline L_i)_{i=0}^d$ be a family of Hermitian
line bundles on $X$. Let $m\in\{0,\cdots d\}$. Assume
that $\overline L_i$ is big for any
$i\in\{0,\cdots,m\}$ and is integrable for any
$i\in\{m+1,\cdots, d\}$. For any $i\in\{0,\cdots,m\}$,
let $D_i=(\nu_i:X_i\rightarrow X,\overline N_i,p_i)$ be
an admissible decomposition of $\overline L_i$. Choose
a birational projective morphism $\nu:X'\rightarrow X$
which factorizes through $\nu_i$ for each
$i\in\{0,\cdots,m\}$. Denote by $\eta_i:X_i\rightarrow
X$ the projective birational morphism such that
$\nu=\nu_i\eta_i$ ($0\leqslant i\leqslant m$). Define
$(D_0\cdots D_{m})\cdot\widehat{c}_1(\overline
L_{m+1})\cdots\widehat{c}_1(\overline L_d)$ as the
normalized intersection product
\[\widehat{c}_1(\eta_0^{*}\overline N_{0})\cdots\widehat{c}_1(
\eta_{m}^*\overline N_{m})\widehat{c}_1(\nu^*\overline
L_{m+1})\cdots\widehat{c}_1(\nu^*\overline
L_d)\prod_{i=0}^mp_i^{-1}.\] This definition does not
depend on the choice of $\nu$.

\begin{prop}\label{Pro:decomposed intersection product}
Let $(\overline{L}_i)_{0\leqslant i\leqslant d}$ be a
family of Hermitian line bundles on $X$. Let
$m\in\{0,\cdots,d\}$. Assume that $\overline{L}_i$ is
big for $i\in\{0,\cdots,m\}$, and is nef for
$i\in\{m+1,\cdots,d\}$. For any $i\in\{0,\cdots,m\}$,
let $D_i$ and $D_i'$ be two admissible decompositions
of $\overline{L}_i$ such that $D_i\succ D_i'$. Then
\begin{equation}\label{Equ:comparison of
product}(D_0\cdots D_{m})\cdot\widehat{c}_1(\overline{
L}_{m+1})\cdots \widehat{c}_1(\overline{
L}_{d})\geqslant (D_0'\cdots
D_{m}')\cdot\widehat{c}_1(\overline{ L}_{m+1})\cdots
\widehat{c}_1(\overline{ L}_{d})\end{equation}
\end{prop}
\begin{proof}
By substituting progressively $D_i$ by $D_i'$, it
suffices to prove that
\[(D_0\cdot D_1\cdots D_{m})\cdot
\widehat{c}_1(\overline{L}_{m+1})\cdots
\widehat{c}_1(\overline{L}_{d})\geqslant (D_0'\cdot
D_1\cdots D_{m})\cdot\widehat{c}_1(\overline{
L}_{m+1})\cdots \widehat{c}_1(\overline{ L}_{d}),\]
which is a consequence of Proposition
\ref{Pro:positivity of nef intersection produc}.
\end{proof}

\begin{coro}
With the notation and the assumptions of Proposition
\ref{Pro:decomposed intersection product}, the supremum
\begin{equation}\label{Equ:intersection product D}\sup\Big\{(D_0\cdots D_{m})
\cdot\widehat{c}_1(\overline{L}_{m+1})\cdots
\widehat{c}_1(\overline{L}_{d})\;\Big|\; 0\leqslant
i\leqslant m,\,
D_i\in\Theta(\overline{L}_i)\Big\}\end{equation} exists
in $\mathbb R_{\geqslant 0}$.
\end{coro}
\begin{proof}
For any $i\in\{0,\cdots,m\}$, let $\overline{A}_i$ be
an arithmetically ample Hermitian line bundle on $X$
such that $\overline{A}_i\otimes\overline{L}_i^{\vee}$
is effective. Then all numbers of the set
\eqref{Equ:intersection product D} is bounded from
above by $\widehat{c}_1(\overline{A}_0)\cdots
\widehat{c}_1(\overline{A}_{m})\cdot\widehat{c}_1(
\overline{L}_{m+1})\cdots\widehat{c}_1(\overline{
L}_{d})$.
\end{proof}

\subsection{Positive intersection product}
\label{SubSec:Posit} Let $(\overline L_i)_{i=0}^{m}$ be
a family of big Hermitian line bundles on $X$, where
$0\leqslant m\leqslant d$. Denote by
$\big\langle\widehat{c}_1(\overline
L_0)\cdots\widehat{c}_1(\overline L_{m})\big\rangle$
the function on $\widehat{\mathrm{Nef}}(X)^{d-m}$ which
sends a family of nef Hermitian line bundles
$(\overline L_j)_{j=m+1}^d$ to the supremum
\[\sup\Big\{(D_0\cdots D_{m})\cdot
\widehat{c}_1(\overline
L_{m+1})\cdots\widehat{c}_1(\overline
L_d)\;\Big|\;0\leqslant i\leqslant m,\,
D_i\in\Theta(\overline L_i) \Big\}.\] Since all
$\Theta(\overline L_i)$ are filtered, this function is
additive in each $\overline L_j$ ($m+1\leqslant
j\leqslant d$). Thus it extends naturally to a
multi-linear function on $\widehat{\mathrm{Int}}(X)$
which we still denote by
$\big\langle\widehat{c}_1(\overline
L_0)\cdots\widehat{c}_1(\overline L_{m})\big\rangle$,
called the {\it positive intersection product} of
$(\overline L_i)_{i=0}^{m}$.

\begin{rema}\label{Rem:property of admi decom}
\begin{enumerate}[(1)]
\item If all Hermitian vector bundles
$(\overline{L}_i)_{i=0}^{m}$ are nef, then the positive
intersection product coincides with the usual
intersection product.
\item The positive intersection product is homogeneous
in each $\overline{L}_i$ ($0\leqslant i\leqslant m$).
However, in general it is not additive in each
variable. If we consider it as a function on
$\widehat{\mathrm{Nef}}(X)$, then it is super-additive
in each variable.
\item Assume that all Hermitian line bundles
$(\overline{L}_i)_{i=0}^{m}$ are the same. That is,
$\overline{L}_0=\cdots=\overline{L}_{m}= \overline{ L
}$. We use the expression
$\big\langle\widehat{c}_1(\overline{L})^{m+1}\big\rangle$
to denote the positive intersection product
\[\big\langle \underbrace{\widehat{c}_1(\overline{
L}) \cdots\widehat{c}_1(\overline{
L})}_{m+1\,\text{copies}}\big\rangle.\] With this
notation, for any $(\overline{
L}_j)_{j=m+1}^{d}\in\widehat{\mathrm{Nef}}( X)^{d-m}$,
one has
\[\begin{split}&\quad\;\big\langle\widehat{c}_1(\overline{L}
)^m\big\rangle\cdot\widehat{c}_1(\overline{
L}_{m+1})\cdots \widehat{c}_1(\overline{
L}_{d})=\sup_{D\in\Theta(\overline{L })}(D\cdots
D)\cdot\widehat{c}_1(\overline{L}_{m+1})\cdots
\widehat{c}_1(\overline{L}_{d}).
\end{split}\]
This equality comes from the fact that the ordered set
$\Theta(\overline{L})$ is filtered (Proposition
\ref{Pro:filtered}) and from the comparison
\eqref{Equ:comparison of product}. In particular, the
Fujita's approximation theorem (see \cite{Chen_Fujita}
and \cite{Yuan08}) implies that
$\big\langle\widehat{c}_1(\overline{ L
})^{d+1}\big\rangle=\widehat{\mathrm{vol}}(\overline{ L
})$.
\end{enumerate}
\end{rema}

\begin{lemm}\label{Lem:comparaison des produVit positive}
Let $(\overline L_i)_{i=0}^m$ be a family of big
Hermitian line bundles on $X$, where
$m\in\{0,\cdots,d\}$. For any $i\in\{0,\cdots,m\}$, let
$\overline M_i$ be an effective Hermitian line bundle
on $X$ and let $\overline N_i=\overline L_i\otimes
\overline M_i$. Then one has
\begin{equation}\label{Equ:comparaison de produit positif}
\big\langle\widehat{c}_1(\overline
L_0)\cdots\widehat{c}_1(\overline
L_m)\big\rangle\geqslant\big\langle\widehat{c}_1(\overline
N_0 )\cdots\widehat{c}_1(\overline N_m)\big\rangle,
\end{equation}
where we have considered the positive intersection
products as functions on
$\widehat{\mathrm{Nef}}(X)^{d-m}$.
\end{lemm}
\begin{proof}
By Remark \ref{Rem:property of admi decom} (4), if
$D_i$ is an admissible decomposition of
$\overline{L}_i$, then it is also an admissible
decomposition of $\overline L_i$. Hence by the
definition of the positive intersection product, the
inequality \eqref{Equ:comparaison de produit positif}
is true.
\end{proof}

The following proposition shows that the positive
intersection product is continuous in each variable.

\begin{prop}\label{Pro:continuite}
Let $(\overline{L}_i)_{0\leqslant i\leqslant m }$ be a
family of big Hermitian line bundles on $X$, where
$m\in\{0,\cdots,d\}$. Let $(\overline{
M}_i)_{0\leqslant i\leqslant m}$ be a family of
Hermitian line bundles on $\mathscr X$. Then
\begin{equation}\label{Equ:continuity}\lim_{n\rightarrow\infty}n^{-m}
\big\langle \widehat{c}_1(\overline{L}_0^{\otimes
n}\otimes\overline{
M}_0)\cdots\widehat{c}_1(\overline{L}_m^{\otimes
n}\otimes\overline{M}_m)\big\rangle=
\big\langle\widehat{c}_1(\overline{L}_0
)\cdots\widehat{c}_1(\overline{L}_m
)\big\rangle\end{equation}
\end{prop}
\begin{proof} We consider firstly both
positive intersection products as functions on
$\widehat{\mathrm{Nef}}(X)$. Let
$\alpha_n=\big\langle\widehat{c}_1(\overline
L_0^{\otimes n} \otimes\overline
M_0)\cdots\widehat{c}_1(\overline L_m^{\otimes n}
\otimes\overline M_m)\big\rangle$. Since $\overline
L_i$ is big, there exists an integer $q\geqslant 1$
such that the Hermitian line bundles $\overline
L_i^{\otimes q}\otimes\overline M_i$ and $\overline
L_i^{\otimes q}\otimes\overline M_i^\vee$ are both
effective. Thus the Lemma \ref{Lem:comparaison des
produVit positive} implies that \begin{gather*}
\alpha_n\geqslant \big\langle\widehat{c}_1(\overline
L_0^{\otimes(n-q)})\cdots\widehat{c}_1(\overline{
L}_m^{\otimes(n-q)})\big\rangle=(n-q)^m\big\langle\widehat{c}_1
(\overline L_0)\cdots\widehat{c}_1(\overline L_m)
\big\rangle,\\
\alpha_n\leqslant \big\langle\widehat{c}_1(\overline
L_0^{\otimes(n+q)})\cdots\widehat{c}_1(\overline{L}_m^{
\otimes(n+q)})\big\rangle=(n+q)^m\big\langle
\widehat{c}_1(\overline
L_0)\cdots\widehat{c}_1(\overline L_m)
\big\rangle.\end{gather*} By passing to limit, we
obtain \eqref{Equ:continuity} as an equality of
functions on $\widehat{\mathrm{Nef}}(X)^{d-m}$. The
general case follows from the multi-linearity.
\end{proof}

\begin{rema}\label{Rem:continuite}
Proposition \ref{Pro:continuite} implies in particular that, if $(f_n^{(i)})_{n\geqslant 1}$ ($i=0,1,\cdots,m$) are families of continuous functions on $X(\mathbb C)$ which converge uniformly to zero. Then one has
\[\lim_{n\rightarrow+\infty}\big\langle\widehat{c}_1(\overline L_0(f_0))\cdots\widehat{c}_m(\overline L_m(f_m))\big\rangle=\big\langle\widehat{c}_1(\overline L_0)\cdots\widehat{c}_m(\overline L_m)\big\rangle.\]
In particular, the mapping
\[t\mapsto\big\langle \widehat{c}_1(\overline L_0(t))\cdots\widehat{c}_1(\overline L_m(t))\big\rangle
\]
is continuous on the (open) interval that it it well defined.
\end{rema}

\begin{prop}\label{Pro:positivite de positive intersection product}
Let $(\overline L_i)_{i=0}^{d-1}$ be a family of big
Hermitian line bundles on $X$. If $\overline M$ is an
effective integrable Hermitian line bundle on $X$, then
\[\big\langle\widehat{c}_1(\overline L_0)\cdots\widehat{c}_1(\overline L_{d-1})\big\rangle
\cdot\widehat{c}_1(\overline M)\geqslant 0.
\]
\end{prop}
\begin{proof}
This is a direct consequence of Proposition
\ref{Pro:positivity of nef intersection produc}.
\end{proof}
\begin{rema}\label{Rem:generations of positive product}
Proposition \ref{Pro:positivite de positive
intersection product} permits us to extend the function
$\big\langle\widehat{c}_1(\overline
L_0)\cdots\widehat{c}_1(\overline L_{d-1})\big\rangle$
on $\widehat{\mathrm{Pic}}(X)$. Let $\overline M$ be an
arbitrary Hermitian line bundle on $X$. By
Weierstrass-Stone theorem, there exists a sequence
$(f_n)_{n\geqslant 1}$ of continuous functions on
$X(\mathbb C)$ which converges uniformly to $0$, and
such that $\overline M(f_n)$ is of smooth metrics for
any $n$. Thus $\overline M(f_n)$ is integrable and
$a_n=\big\langle\widehat{c}_1(\overline
L_0)\cdots\widehat{c}_1(\overline
L_{d-1})\big\rangle\cdot\widehat{c}_1(\overline
M(f_n))$ is well defined. Let
$\varepsilon_{n,m}=\|f_n-f_m\|_{\sup}$. Choose an ample
Hermitian line bundle $\overline A$ such that
$\overline A\otimes\overline L_{i}^\vee$ is effective
for any $i\in\{0,\cdots,d-1\}$. Note that
\begin{equation}\label{Equ:difference a m et an}\begin{split}a_n-a_m&=\big\langle\widehat{c}_1(\overline
L_0)\cdots\widehat{c}_1(\overline
L_{d-1})\big\rangle\cdot\widehat{c}_1(\overline{\mathcal
O }(f_n-f_m))\\
&\leqslant \langle\widehat{c}_1(\overline
L_0)\cdots\widehat{c}_1(\overline
L_{d-1})\big\rangle\cdot\widehat{c}_1(\overline{\mathcal
O }(\varepsilon_{n,m}))\\
&\leqslant\big\langle\widehat{c}_1(\overline
A)^d\big\rangle\cdot\widehat{c}_1(\overline{\mathcal
O}(\varepsilon_{n,m}))=\varepsilon_{n,m}c_1(A_K)^d
\end{split}\end{equation}
By interchanging the roles of $n$ and $m$ in
\eqref{Equ:difference a m et an} and then combining the
two inequalities, one obtains
$|a_n-a_m|\leqslant\varepsilon_{n,m}c_1(A_K)^d$.
Therefore, $(a_n)_{n\geqslant 1}$ is a Cauchy sequence
which converges to a real number which we denote by
$\big\langle\widehat{c}_1(\overline
L_0)\cdots\widehat{c}_1(\overline
L_{d-1})\big\rangle\cdot\widehat{c}_1(\overline M)$. By
an argument similar to the inequality
\eqref{Equ:difference a m et an}, this definition does
not depend on the choice of the sequence
$(f_n)_{n\geqslant 1}$. The extended function is
additive on $\widehat{\mathrm{Pic}}(X)$, which is
positive on the subgroup of effective Hermitian line
bundles, and satisfies the conclusion of Proposition \ref{Pro:continuite}.
\end{rema}

\section{Differentiability of the arithmetic volume function}

In this section, we establish the differentiability of
the arithmetic volume function.
We begin by presenting the following lemma, which is
analogous to \cite[Corollary 3.4]{Bou_Fav_Mat06}.
\begin{lemm}\label{Lem:Siu-Yuan}
Let $\overline L$ and $\overline N$ be two nef
Hermitian line bundles on $X$. Let $\overline M $ be an
integrable Hermitian line bundle on $X$. Assume that
$\overline M\otimes\overline N $ and $\overline
M^\vee\otimes\overline N$ are nef and that $\overline
L^\vee\otimes\overline N$ is effective. Then there
exists a constant $C>0$ only depending on $d$ such that
\[\widehat{\mathrm{vol}}(\overline L^{\otimes n}
\otimes\overline M)\geqslant
n^{d+1}\widehat{\mathrm{vol}}(\overline L)
+(d+1)n^d\widehat{c}_1(\overline
L)^{d}\widehat{c}_1(\overline
M)-C\widehat{\mathrm{vol}}(\overline N)n^{d-1}.\]
\end{lemm}
\begin{proof}
Recall that (see \cite[Theorem 2.2]{Yuan07}, see also \cite[Theorem 5.6]{Moriwaki07}) if
$\overline A$ and $\overline B$ are two nef Hermitian
line bundles on $X$, then
\begin{equation}\label{Equ:Yuan's inequality}\widehat{\mathrm{vol}}(\overline B\otimes\overline
A^\vee) \geqslant \widehat{c}_1(\overline
B)^{d+1}-\widehat{c}_1(\overline
B)^d\widehat{c}_1(\overline A).\end{equation} Let
$\overline B=\overline L^{\otimes n}\otimes\overline
M\otimes\overline N$. It is a nef Hermitian line bundle
on $X$. If one applies \eqref{Equ:Yuan's inequality} on
$\overline B$ and on $\overline A=\overline N$, one
obtains
\[\begin{split}\widehat{\mathrm{vol}}(\overline L^{\otimes n}
\otimes\overline M)&=\widehat{\mathrm{vol}} (\overline
B\otimes\overline N^\vee)\geqslant
\widehat{c}_1(\overline
B)^{d+1}-(d+1)\widehat{c}_1(\overline
B)^d\widehat{c}_1(\overline N)\\
&=n^{d+1}\widehat{c}_1(\overline L)^{d+1}+
(d+1)n^d\widehat{c}_1(\overline
L)^d\widehat{c}_1(\overline M)+O(n^{d-1}),
\end{split}
\]
where the implicit constant is a linear combination of
intersection numbers of Hermitian line bundles of the
form $\overline L $ or $\overline M\otimes\overline N$,
and hence can be bounded from above by a multiple of
$\widehat{c}_1(\overline
N)^{d+1}=\widehat{\mathrm{vol}}(\overline N)$,
according to Proposition \ref{Pro:positivity of nef
intersection produc}.
\end{proof}

\begin{rema}
In \cite{Yuan07}, Yuan has actually proved a stronger
inequality by replacing the
$\widehat{\mathrm{vol}}(\overline B\otimes\overline
A^\vee)$ in \eqref{Equ:Yuan's inequality} by
\[\lim_{n\rightarrow+\infty}\frac{\chi(\pi_*(\overline
B^{\otimes n}\otimes\overline A^{\vee\otimes
n}))}{n^{d+1}/(d+1)!}.\] In fact, this quantity is
always bounded from above by
$\widehat{\mathrm{vol}}(\overline B\otimes\overline
A^\vee)$.
\end{rema}

\begin{proof}[Proof of Theorem \ref{Thm:main theorem}] We first assume that the metrics of $\overline M$ are
smooth. Choose an ample Hermitian line bundle
$\overline N$ such that $\overline N\otimes\overline M$
and $\overline N\otimes\overline M^\vee$ are ample, and
that $\overline N\otimes\overline L^\vee $ is
effective. Let $D=(\nu:X'\rightarrow X,\overline A,p)$
be an admissible decomposition of $\overline L$. One
has
\[\widehat{\mathrm{vol}}(\overline L^{\otimes n}
\otimes\overline M)=p^{-d-1}
\widehat{\mathrm{vol}}(\nu^*\overline L^{\otimes
np}\otimes\nu^*\overline M^{\otimes p})\geqslant
p^{-d-1}\widehat{\mathrm{vol}}(\overline A^{\otimes
n}\otimes\overline{M}^{\otimes p }).\] Note that
$\overline N^{\otimes p}\otimes\overline M^{\otimes p
}$ and $\overline N^{\otimes p}\otimes\overline
M^{\vee\otimes p}$ are nef and \[\overline N^{\otimes p
}\otimes\overline A^{\vee}=(\overline N\otimes
\overline L^{\vee})^{\otimes p}\otimes(\overline
L^{\otimes p }\otimes\overline A^\vee)\] is effective.
After Lemma \ref{Lem:Siu-Yuan}, one obtains
\[\widehat{\mathrm{vol}}(\overline A^{\otimes n}
\otimes\overline M^{\otimes p}) \geqslant
n^{d+1}\widehat{c}_1(\overline A)^{d+1}+
(d+1)n^dp\widehat{c}_1(\overline A)^d
\widehat{c}_1(\overline
M)-Cp^{d+1}\widehat{\mathrm{vol}}(\overline N
)n^{d-1}.\] Therefore,
\[\widehat{\mathrm{vol}}(\overline L^{\otimes n}
\otimes\overline M)\geqslant
n^{d+1}(D^{d+1})+(d+1)n^d(D^d)\cdot\widehat{c}_1(\overline
M)-C\widehat{\mathrm{vol}}(\overline N )n^{d-1}.\]
Since $D $ is arbitrary, one has
\begin{equation}\label{Equ:estimation de siu} \widehat{\mathrm{vol}}(\overline
L^{\otimes n} \otimes\overline M)\geqslant
n^{d+1}\widehat{\mathrm{vol}}(\overline L
)+(d+1)n^d\big\langle\widehat{c}_1(\overline L
)^d\big\rangle\cdot\widehat{c}_1(\overline
M)-C\widehat{\mathrm{vol}}(\overline N)n^{d-1}.
\end{equation}
By passing to limit, one obtains
\[\liminf_{n\rightarrow+\infty}\frac{\widehat{\mathrm{vol}}
(\overline L^{\otimes n}\otimes\overline
M)-\widehat{\mathrm{vol}}(\overline L^{\otimes
n})}{n^d}\geqslant
(d+1)\big\langle\widehat{c}_1(\overline L
)^d\rangle\cdot\widehat{c}_1(\overline M).\] If we
apply \eqref{Equ:estimation de siu} in replacing
$\overline L$ by $\overline L^{\otimes
n}\otimes\overline M$, $\overline M$ by $\overline
M^{\vee\otimes n}$ and $\overline N$ by $\overline
N^{\otimes 2n}$, we obtain
\[\widehat{\mathrm{vol}}(\overline L^{\otimes n^2})
\geqslant n^{d+1}\widehat{\mathrm{vol}}(\overline
L^{\otimes n}\otimes \overline
M)-(d+1)n^{d+1}\big\langle\widehat{c}_1(\overline
L^{\otimes n}\otimes\overline
M)^d\big\rangle\cdot\widehat{c}_1(\overline
M)-C(2n)^{d+1}\widehat{\mathrm{vol}}(\overline N)
n^{d-1}, \] or equivalently
\[\widehat{\mathrm{vol}}(\overline L^{\otimes n})
\geqslant\widehat{\mathrm{vol}}(\overline L^{\otimes
n}\otimes\overline M)-
(d+1)\big\langle\widehat{c}_1(\overline L^{\otimes
n}\otimes\overline
M)^d\big\rangle\cdot\widehat{c}_1(\overline
M)-2^{d+1}C\widehat{\mathrm{vol}}(\overline
N)n^{d-1}.\] Thus
\[\begin{split}\limsup_{n\rightarrow\infty}\frac{\widehat{\mathrm{vol}}
(\overline L^{\otimes n}\otimes\overline
M)-\widehat{\mathrm{vol}}(\overline L^{\otimes
n})}{n^d}&\leqslant\lim_{n\rightarrow\infty}n^{-d}(d+1)
\big\langle\widehat{c}_1(\overline L^{\otimes n}\otimes\overline M)^d\big\rangle\cdot\widehat{c}_1(\overline M)\\
&=(d+1)\big\langle\widehat{c}_1(\overline
L)^d\big\rangle\cdot\widehat{c}_1(\overline M).
\end{split}\]
Therefore, one has
\[\lim_{n\rightarrow+\infty}\frac{\widehat{\mathrm{vol}}
(\overline L^{\otimes n}\otimes\overline
M)-\widehat{\mathrm{vol}}(\overline L^{\otimes
n})}{n^d}=(d+1)\big\langle\widehat{c}_1(\overline L
)^d\rangle\cdot\widehat{c}_1(\overline M).\]

For the general case, by Weierstrass-Stone theorem, for
any $\varepsilon>0$, there exist two Hermitian line
bundles with smooth metrics $\overline
M_{\varepsilon,1}=(M,(\|\cdot\|_{\sigma,\varepsilon}')_{\sigma:K\rightarrow\mathbb
C })$ and $\overline
M_{\varepsilon,2}=(M,(\|\cdot\|_{\sigma,\varepsilon}'')_{\sigma:K\rightarrow\mathbb
C })$ such that
\[\|\cdot\|_{\sigma,\varepsilon}'\leqslant\|\cdot\|_\sigma
\leqslant\|\cdot\|_{\sigma,\varepsilon}'',\quad\text{and}\quad
\max_{\sigma}\sup_{x\in X_\sigma(\mathbb C)} \big|
\log\|\cdot\|_{\sigma,\varepsilon}'(x)-
\log\|\cdot\|_{\sigma,\varepsilon}''(x)\big|
\leqslant\varepsilon,\] where $\|\cdot\|_\sigma$ is the
norm of index $\sigma$ of $\overline M$. Note that
$\overline M_{\varepsilon,2}\leqslant\overline
M\leqslant\overline M_{\varepsilon,1} $. By the special
case that we have proved, one has
\[\begin{split}&\quad\;(d+1)\big\langle\widehat{c}_1(\overline L)^d\big\rangle\cdot
\widehat{c}_1(\overline
M_{\varepsilon,2})\leqslant\liminf_{n\rightarrow+\infty}\frac{\widehat{\mathrm{vol}}
(\overline L^{\otimes n}\otimes\overline
M)-\widehat{\mathrm{vol}}(\overline L^{\otimes
n})}{n^d}\\&\leqslant\limsup_{n\rightarrow+\infty}\frac{\widehat{\mathrm{vol}}
(\overline L^{\otimes n}\otimes\overline
M)-\widehat{\mathrm{vol}}(\overline L^{\otimes
n})}{n^d}\leqslant
(d+1)\big\langle\widehat{c}_1(\overline L
)^d\big\rangle\cdot\widehat{c}_1(\overline
M_{\varepsilon,1} )\end{split}\] Let
$f_\varepsilon:X(\mathbb C)\rightarrow\mathbb R$ be the
function such that $\log\|\cdot\|_{\sigma,\varepsilon}''(x)=
\log\|\cdot\|_{\sigma,\varepsilon}'(x)+f_\varepsilon(x)$.
Denote by $\overline{\mathcal O}(f_\varepsilon)$ the
Hermitian line bundle on $X$ whose underlying Hermitian
line bundle is trivial, and such that
$\|\mathbf{1}\|(x)=e^{-f_\varepsilon(x)}$. It is an
effective Hermitian line bundle since
$f_{\varepsilon}\geqslant 0$. Furthermore, one has
$\overline M_{\varepsilon,1}\otimes\overline
M_{\varepsilon, 2}^\vee\cong\overline{\mathcal
O}(f_\varepsilon)$. Let $\overline F$ be an ample
Hermitian line bundle on $X$ such that $\overline
F\otimes\overline L^\vee$ is effective. One has
\[\begin{split}&\quad\;\big\langle\widehat{c}_1(\overline L)^d\big\rangle
\cdot(\widehat{c}_{1}(\overline
M_{\varepsilon,2})-\widehat{c}_1(\overline
M_{\varepsilon,1}))=\big\langle\widehat{c}_1(\overline
L)^d\big\rangle\cdot\widehat{c}_1(\overline {\mathcal
O}(f_\varepsilon) )\\
&\leqslant\big\langle\widehat{c}_1(\overline
L)^d\big\rangle\cdot\widehat{c}_1(\overline {\mathcal
O}(\varepsilon) )\leqslant \widehat{c}_1(\overline
F)^d\widehat{c}_1(\overline {\mathcal O}(\varepsilon)
)= \varepsilon\int_{X(\mathbb C)}{c}_1(\overline
F)^d,\end{split}\] where in the first inequality, we
have applied Proposition \ref{Pro:positivite de
positive intersection product} (see also Remark
\ref{Rem:generations of positive product}), and in the
second inequality, we have used the fact that
$\overline{\mathcal O}(\varepsilon)$ is nef and then
applied Lemme \ref{Lem:comparaison des produVit
positive}. Since $\varepsilon$ is arbitrary, we obtain
that
\[ D_{\overline L }\widehat{\mathrm{vol}}(\overline M
):=\lim_{n\rightarrow+\infty}\frac{\widehat{\mathrm{vol}}
(\overline L^{\otimes n}\otimes\overline
M)-\widehat{\mathrm{vol}}(\overline L^{\otimes
n})}{n^d}\] exists in $\mathbb R$. Since
\[(d+1)\big\langle\widehat{c}_1(\overline L)^d\big\rangle\cdot
\widehat{c}_1(\overline M_{\varepsilon,2})\leqslant
D_{\overline L }\widehat{\mathrm{vol}}(\overline M
)\leqslant (d+1)\big\langle\widehat{c}_1(\overline L
)^d\big\rangle\cdot\widehat{c}_1(\overline
M_{\varepsilon,1} ),\] by virtue of Remark
\ref{Rem:generations of positive product}, we obtain $
D_{\overline L }\widehat{\mathrm{vol}}(\overline M
)=(d+1)\big\langle\widehat{c}_1(\overline L
)^d\big\rangle\cdot\widehat{c}_1(\overline M)$, which
is additive with respect to $\overline M$.
\end{proof}

A direct consequence of Theorem \ref{Thm:main theorem}
is the asymptotic orthogonality of arithmetic Fujita
approximation.
\begin{coro} Assume that $\overline L$ is a big
Hermitian line bundle on $X$. One has
\begin{equation}\label{Equ:vol comm intersection}\big\langle \widehat{c}_1(\overline
L)^d\big\rangle\cdot\widehat{c}_1(\overline L)
=\widehat{\mathrm{vol}}(\overline L).\end{equation}
\end{coro}
\begin{proof}
By definition, \[
\begin{split}
&\quad\;D_{\overline L}\widehat{\mathrm{vol}}(\overline
L)=\lim_{n\rightarrow+\infty}\frac{\widehat{\mathrm{vol}}(\overline
L^{n+1} )-\widehat{\mathrm{vol}}(\overline L^{\otimes n
})}{n^d}\\&=\widehat{\mathrm{vol}}(\overline L)
\lim_{n\rightarrow+\infty}\frac{(n+1)^{d+1}-n^{d+1}}{n^d}
=(d+1)\widehat{\mathrm{vol}}(\overline L).
\end{split}\]
So \eqref{Equ:vol comm intersection} follows from
Theorem \ref{Thm:main theorem}.
\end{proof}

\begin{rema}
As mentioned in Introduction, the differentiability of the geometrical volume function can be obtained by using the method of Okounkov bodies developed in \cite{Okounkov96}. See \cite{Lazarsfeld_Mustata08} for a proof of this result and other interesting results concerning geometric volume functions. Recently, Yuan \cite{Yuan08} has proposed a partial analogue of the construction of Lazarsfeld and Mus\c{t}at\v{a} in Arakelov geometry. In fact, he has defined the Okounkov bodies of a Hermitian line bundle with respect to so-called ``vertical flags''. This permits him to obtain the log-concavity of the arithmetic volume function and the analogue of Fujita's approximation theorem in Arakelov goemetry, where the latter has also been independently  obtained by the present author \cite{Chen_Fujita}, using asymptotic measures. Quite possibly, an analogue of Lazarsfeld and Mus\c{t}at\v{a}'s construction with respect to ``horizontal flags'' could also imply the differentiability of the arithmetic volume function.
\end{rema}

\section{Applications and comparisons}
In this section, we shall apply our differentiability
result to study several arithmetic invariants of
Hermitian line bundles.
\subsection{Asymptotic measure}
Let $\overline L$ be a Hermitian line bundle on $X$
such that $L_K$ is big. The {\it asymptotic measure} of
$\overline L$ is  the vague limit  in the space of
Borel probability measures
\begin{equation}\label{Equ:nu}\nu_{\overline L}:=-\lim_{n\rightarrow+\infty}\frac{\mathrm
d}{\mathrm{d}t}
\frac{\rang\Big(\mathrm{Vect}_K\big(\{s\in\pi_*L^{\otimes
n }\mid\forall\sigma,\,\|s\|_{\sigma,\sup}\leqslant
e^{-\lambda n} \}\big)\Big)}{\rang(\pi_*L^{\otimes n
})},\end{equation} where the derivative is taken in the
sense of distribution. It is also the limit of
normalized Harder-Narasimhan measures (cf.
\cite{Chen08,Chen_bigness,Chen_Fujita}).

Note that the support of the probability measure $\nu_{\overline L}$
is contained in
$]-\infty,\widehat{\mu}_{\max}^\pi(\overline L)]$,
where $\widehat{\mu}_{\max}^\pi(\overline L)$ is the
limit of maximal slopes (see \cite[Theorem
4.1.8]{Chen08}):
\[\widehat{\mu}_{\max}^{\pi}(\overline L):=\lim_{n\rightarrow+\infty}
\frac{\widehat{\mu}_{\max}(\pi_*\overline L^{\otimes
n})}{n}.\] Recall that in \cite[Theorem
5.5]{Chen_bigness}, the present author has proved that
$\widehat{\mu}_{\max}^{\pi}(\overline L)>0$ if and only
if $L_K$ is big. Furthermore, by definition, one has
$\widehat{\mu}_{\max}^{\pi}(\overline
L(a))=\widehat{\mu}_{\max}^{\pi}(\overline L)+a$ for
any $a\in\mathbb R$. Therefore,
$\widehat{\mu}_{\max}^{\pi}(\overline L)$ is also the
infimum of all real numbers $\varepsilon$ such that
$\overline L(-\varepsilon)$ is big.

The asymptotic measure is a very general arithmetic
invariant. Many arithmetic invariants of $\overline L$
can be represented as integrals with respect to
$\nu_{\overline L}$. In the following, we discuss some
examples. The {\it asymptotic positive slope} of
$\overline L$ is defined as
\[\widehat{\mu}_+^\pi(\overline L):=\frac{1}{[K:\mathbb Q]}
\frac{\widehat{\mathrm{vol}}(\overline
L)}{(d+1)\mathrm{vol}(L_K)}.\] In \cite{Chen_bigness},
the author has proved that
$\widehat{\mu}_+^\pi(\overline L)$ is also the maximal
value of the asymptotic Harder-Narasimhan polygon of
$\overline L$ and that the asymptotic positive slope
has the following integral form:
\begin{equation}
\widehat{\mu}_+^\pi(\overline L)=\int_{\mathbb
R}\max(x,0)\,\nu_{\overline L}(\mathrm{d}x).
\end{equation}
More generally, for any $a\in\mathbb R$, one has
\begin{equation}\label{Equ:call et mu plus}
\int_{\mathbb R}\max(x-a,0)\,\nu_{\overline
L}(\mathrm{d}x)=\widehat{\mu}_+^\pi(\overline L(-a)).
\end{equation}

Another important example is the {\it asymptotic slope}
of $\overline L$, which is
\[\widehat{\mu}^\pi(\overline L):=\frac{1}{[K:\mathbb Q]}
\frac{S(\overline
L)}{(d+1)\mathrm{vol}(L_K)}\in[-\infty,+\infty[,\]
where $S(\overline L)$ is the sectional capacity of
$\overline L$ as in \eqref{Equ:sectional
capa}. The asymptotic slope has the following integral
form \begin{equation*}\widehat{\mu}^\pi(\overline
L)=\int_{\mathbb R}x\,\nu_{\overline L}
(\mathrm{d}x).\end{equation*} Observe that we have
\begin{equation}\widehat{\mu}_{\max}^\pi(\overline L)\geqslant
\widehat{\mu}_+^\pi(\overline
L)\geqslant\widehat{\mu}^\pi(\overline
L).\end{equation}

Using Theorem \ref{Thm:main theorem} and the
differentiability of geometric volume function in
\cite{Bou_Fav_Mat06}, we prove that the asymptotic
positive slope $\widehat{\mu}_+^\pi$ is differentiable
and calculate its differential.
\begin{prop}\label{Pro:derivabilite de mu plus}
Assume that $\overline L$ is a big Hermitian line
bundle on $X$. For any Hermitian line bundle $\overline
M$, one has
\[D_{\overline L}\widehat{\mu}_+^\pi(\overline M):=\lim_{n\rightarrow+\infty}
\big({ \widehat{\mu}_+^\pi(\overline L^{\otimes
n}\otimes\overline M )-\widehat{\mu}_+^\pi(\overline
L^{\otimes n})}\big)\] exists in $\mathbb R$.
Furthermore, one has
\[D_{\overline L}\widehat{\mu}_+^\pi(\overline M):=
\frac{\big\langle\widehat{c}_1(\overline
L)^d\big\rangle\cdot\widehat{c}_1(\overline M
)}{[K:\mathbb
Q]\mathrm{vol}(L_K)}-\frac{d\big\langle
c_1(L_K)^{d-1}\big\rangle\cdot c_1(M_K) }{\mathrm{vol}(L_K)}\widehat{\mu}_+^\pi(\overline L),\] where $\big\langle
c_1(L_K)^{d-1}\big\rangle\cdot c_1(M_K)$ is the
geometric positive intersection product (\cite[\S
2]{Bou_Fav_Mat06}).
\end{prop}
\begin{proof}
This is a direct consequence of Theorem \ref{Thm:main
theorem} and \cite[Theorem A]{Bou_Fav_Mat06}, where the
latter asserts that
\[\lim_{n\rightarrow+\infty}\frac{\mathrm{vol}(L_K^{\otimes n}
\otimes M_K)-\mathrm{vol}(L_K^{\otimes n})}{n^{d-1}}=d\big\langle
c_1(L_K)^{d-1}\big\rangle\cdot c_1(M_K).\]
\end{proof}

We then deduce from Proposition \ref{Pro:derivabilite de mu plus} the expression of the distribution function of the measure $\nu_{\overline L}$.
\begin{prop}\label{Pro:distribution function}
The distribution function $F_{\overline L}$ of
$\nu_{\overline L}$ satisfies the equality
\[F_{\overline L}(a):=\nu_{\overline L}(]-\infty,a])=1-\frac{\big\langle\widehat{c}_1(\overline L(-a))^d\big\rangle\cdot\widehat{c}_1(\overline{\mathcal O}(1))}{[K:\mathbb Q]\mathrm{vol}(L_K)}, \quad a<\widehat{\mu}_{\max}^\pi(\overline L).
\]
\end{prop}
\begin{proof}
One has
\[F_{\overline L}(a)=1+\frac{\mathrm{d}}{\mathrm{d}a}\int_{\mathbb R}\max(x-a,0)\,{\nu}(\mathrm{d}x)=1+\frac{\mathrm{d}}{\mathrm{d}a}
\widehat{\mu}_+^\pi(\overline L(-a)).
\]
By Proposition \ref{Pro:derivabilite de mu plus}, one obtains
\[\frac{\mathrm{d}}{\mathrm{d}a}\widehat{\mu}_+^\pi(\overline L(-a))=-\frac{\big\langle\widehat{c}_1(\overline L(-a))^d\big\rangle\cdot\widehat{c}_1(\overline{\mathcal O}(1))}{[K:\mathbb Q]\mathrm{vol}(L_K)}.\]
\end{proof}

\begin{rema}\begin{enumerate}[1)]
\item Since the support of $\nu_{\overline L}$ is
bounded from above by $\widehat{\mu}_{\max}(\overline L
)$, one has $F_{\overline L}(a)=1$ for
$a\geqslant\widehat{\mu}_{\max}(\overline L )$.
\item As a consequence of Proposition
\ref{Pro:distribution function}, one obtains that the
function
\[\frac{\big\langle\widehat{c}_1(\overline L(-a))^d\big\rangle\cdot\widehat{c}_1(\overline{\mathcal O}(1))}{[K:\mathbb Q]\mathrm{vol}(L_K)}\]
is decreasing with respect to $a$ on
$]-\infty,\widehat{\mu}_{\max}^\pi(\overline L)[$,
which is also implied by Lemma \ref{Lem:comparaison des
produVit positive}. Furthermore, this function takes
values in $]0,1]$, and converges to $1$ when $a\rightarrow-\infty$.
\item Let $a\in]-\infty,\widehat{\mu}_{\max}^\pi(\overline L)[$. The restriction of
\[\frac{1}{[K:\mathbb Q]\mathrm{vol}(L_K)}\big\langle
\widehat{c}_1(\overline L(-a))^d\big\rangle\] on
$C^0(X(\mathbb C))$ (considered as a subgroup of
$\widehat{\Pic}(X)$ via the mapping $f\mapsto\mathcal
O(f)$) is a positive linear functional, thus
corresponds to a Radon measure on $X(\mathbb C)$.
Furthermore, by 1), its total mass is bounded from
above by $1$, and converges to $1$ when $a\rightarrow-\infty$.
\item After Remark
\ref{Rem:continuite}, we observe from Proposition
\ref{Pro:distribution function} that the only possible
discontinuous point of the distribution function
$F_{\overline L}(a)$ is
$a=\widehat{\mu}_{\max}^\pi(\overline L)$.
\end{enumerate}
\end{rema}

As an application, we calculate the sectional capacity in terms of positive intersection product.
\begin{coro} Let $\overline L$ be a Hermitian line bundle on $X$ such that $L_K$ is big. Let $A=\widehat{\mu}_{\max}^\pi(\overline L)$.
One has
\[S(\overline L)=(d+1)A\lim_{x\rightarrow A-}\big\langle\widehat{c}_1(\overline L(-x))^d\big\rangle\cdot\widehat{c}_1(\overline {\mathcal O}(1))-\int_{-\infty}^{A}
(d+1)x\,\mathrm{d}\big\langle\widehat{c}_1(\overline L(-x))^d\big\rangle\cdot\widehat{c}_1(\overline{\mathcal O}(1))
.\]
\end{coro}

\subsection{Lower bound of the positive intersection
product} Our differentiability result permits to obtain
a lower bound for positive intersection products of the
form $\big\langle\widehat{c}_1(\overline
L)^d\big\rangle\cdot\widehat{c}_1(\overline M)$, where
$\overline L$ is a big Hermitian line bundle on $X$ and
$\overline M$ is an effective Hermitian lien bundle on
$X$, by using the log-concavity of the arithmetic
volume function proved in \cite{Yuan08}.

\begin{prop}
Let $\overline L$ and $\overline M$ be two Hermitian line bundles on $X$. Assume that $\overline L$ is big and $\overline M$ is effective. Then
\begin{equation}\label{Equ:isoperimetric}\big\langle\widehat{c}_1(\overline L)^d\big\rangle\cdot\widehat{c}_1(\overline M)\geqslant
\mathrm{vol}(\overline L)^{\frac{d}{d+1}}\mathrm{vol}(\overline M)^{\frac 1{d+1}}.\end{equation}
\end{prop}
\begin{proof}
Theorem \ref{Thm:main theorem} shows
\[\lim_{n\rightarrow+\infty}\frac{\mathrm{vol}(\overline L^{\otimes n}\otimes\overline M)-\mathrm{vol}(\overline M)}{n^d}=(d+1)\big\langle\widehat{c}_1(\overline L)^d\big\rangle\cdot\widehat{c}_1(\overline M).\]
By \cite[Theorem B]{Yuan08}, one has
\[\widehat{\mathrm{vol}}(\overline L^{\otimes n}\otimes\overline M)\geqslant\Big(\widehat{\mathrm{vol}}(\overline L^{\otimes n})^{\frac{1}{d+1}}+\widehat{\mathrm{vol}}(\overline M)^{\frac{1}{d+1}}\Big)^{d+1}.\]
By passing to limit, we obtain the required inequality.
\end{proof}

\begin{rema}
The inequality \eqref{Equ:isoperimetric} could be considered as an analogue in Arakelov geometry (suggested by Bertrand \cite{Bertrand95}) of the {\it isoperimetric inequality} proved by Federer \cite[3.2.43]{Federer69}. See \cite[\S 5.4]{Fulton93} for an interpretation in terms of intersection theory, and \cite[\S 1.2]{Bertrand95} for an analogue in geometry of numbers.
\end{rema}

\subsection{Comparison to other differentiability
results}

We finally compare our results to several
differentiability results on arithmetic invariants.

\subsubsection*{Intersection number} Recall that the
self-intersection number $\widehat{c}_1(\overline
L)^{d+1}$ is well defined for integrable Hermitian line
bundles $\overline L$. See
\cite{Gillet-Soule90,Zhang95,Zhang95b}. Furthermore, it
is a polynomial function. Therefore, for any integrable
Hermitian line bundles $\overline L$ and $\overline M$,
one has
\[\lim_{n\rightarrow+\infty}\frac{\widehat{c}_1(\overline L^{\otimes n}\otimes\overline M)^{d+1}-\widehat{c}_1(\overline L^{\otimes n})}{n^d}
=(d+1)\widehat{c}_1(\overline
L)^{d}\widehat{c}_1(\overline M).\] This formula shows
that the intersection number is differentiable at
$\overline L$ along all directions in
$\widehat{\mathrm{Int}}(X)$.

\subsubsection*{Sectional capacity}

By using the analogue of Siu's inequality in Arakelov geometry, Yuan \cite{Yuan07} has actually proved that the sectional capacity $S$ is differentiable  along integrable directions
 at any Hermitian line bundle $\overline L$ such that $L$ is ample and that the metrics of $\overline L$ are semi-positive. Furthermore, for such $\overline L$, one has
\[D_{\overline L}S(\overline M):=\lim_{n\rightarrow+\infty}
\frac{S(\overline L^{\otimes n}\otimes\overline M)-S(\overline L^{\otimes n})}{n^{d}}=(d+1)\widehat{c}_1(\overline L)^d\widehat{c}_1(\overline M),
\]
where $\overline M$ is an arbitrary integrable
Hermitian line bundles. This result has been
established by Autissier \cite{Autissier01} in the case
where $d=1$. Recently Berman and Boucksom
\cite{Ber_Bou08} have proved a general
differentiability result for the sectional capacity.
They have proved that the function $S$ is
differentiable along the directions defined by
continuous functions on $X(\mathbb C)$ on the cone of
generically big Hermitian lien bundles. Namely, for any
continuous function $f$ on $X(\mathbb C)$ and any
Hermitian line bundle $\overline L$ on $X$ such that
$L_K$ is big and that $S(\overline L)$ is finite, the
limit
\[\lim_{n\rightarrow+\infty}\frac{S(\overline L^{\otimes n}(f))-S(\overline L^{\otimes n})}{n^d}\]
exists. They have also computed explicitly the
differential in terms of the Monge-Amp\`ere measure of
$\overline L$ (see Theorem 5.7 and Remark 5.8 {\it loc.
cit.}).

Our differentiability result for arithmetic volume
function (Theorem \ref{Thm:main theorem}), combined
with \eqref{Equ:call et mu plus}, implies the
differentiability of arithmetic invariants which can be
written as the integration of a (fixed) smooth function
of compact support with respect to the asymptotic
measure of the Hermitian line bundle, by using
integration by part. It would be interesting to know if
a similar idea permits to deduce the differentiability
of the sectional capacity, which can be written as the
integral of the function $f(x)=x$ with respect to the
asymptotic measure, along any direction at any
Hermitian line bundle $\overline L$ such that $L_K$ is
big and that $S(\overline L)$ is finite.

\backmatter
\bibliography{chen}
\bibliographystyle{smfplain}

\end{document}